# GLOBALLY OPTIMAL PARAMETER ESTIMATES FOR NONLINEAR DIFFUSIONS

By Aleksandar Mijatović and Paul Schneider

*Imperial College London and Warwick Business School*

This paper studies an approximation method for the log-likelihood function of a nonlinear diffusion process using the bridge of the diffusion. The main result (Theorem 1) shows that this approximation converges uniformly to the unknown likelihood function and can therefore be used efficiently with any algorithm for sampling from the law of the bridge. We also introduce an expected maximum likelihood (EML) algorithm for inferring the parameters of discretely observed diffusion processes. The approach is applicable to a subclass of nonlinear SDEs with constant volatility and drift that is linear in the model parameters. In this setting, globally optimal parameters are obtained in a single step by solving a linear system. Simulation studies to test the EML algorithm show that it performs well when compared with algorithms based on the exact maximum likelihood as well as closed-form likelihood expansions.

**1. Introduction.** In the natural and social sciences, diffusion processes are widely used as models for random phenomena that evolve continuously in time. They are popular because they arise as solutions of stochastic differential equations, which are natural probabilistic generalizations of the deterministic models described by ordinary differential equations. It is well known that if the data are recorded at discrete times, parametric inference for diffusions using maximum likelihood estimates is difficult, primarily because it is usually impossible to find the corresponding likelihood function [see Sørensen (2004) for a review of methods of inference in the diffusion setting]. In this paper, we are concerned with the estimation of the parameters in the stochastic differential equation (SDE)

$$(1.1) \qquad dX_t = \mu(X_t, \theta)\,dt + dW_t \qquad \text{where } \mu(\cdot, \cdot)\colon \mathbb{R} \times \mathbb{R}^{N+1} \to \mathbb{R},$$

is an arbitrary continuous (possibly) nonlinear function and $W_t$ is a standard Brownian motion. In order to guarantee the existence of a nonexplosive solution of (1.1), we need to assume that for each parameter value $\theta \in \mathbb{R}^{N+1}$, the function $\mu(\cdot, \theta) : \mathbb{R} \to \mathbb{R}$ is locally Lipschitz with linear growth [see Kloeden and Platen (1999), Chapter 4]. The task is to infer the vector of coefficients $\theta \in \Theta \subseteq \mathbb{R}^{N+1}$ in the drift $\mu(\cdot, \theta)$ from $K+1$ observed realizations $x_0, \ldots, x_K$ of the diffusion $X_t$, where $\Theta$ is some compact subset in the parameter space.

When the exact likelihood function is available, the parameters can be determined by maximizing the joint likelihood of the observations. The true likelihood function is, however, available only in very few cases. A variety of approximations exist and are well documented in the literature [see Aït-Sahalia (2002) and Jensen and Poulsen (2002), Hurn, Jeisman and Lindsay (2007), Schneider (2006) for empirical comparisons of the available methods]. A general method for parameter inference based on the EM algorithm is to maximize an approximate likelihood function, which can be defined if one can simulate the bridge of the diffusion in (1.1) (see Section 2 for the precise definition of this approximation). Recently, an exact simulation approach for diffusion bridges was developed in Beskos et al. (2006) and an efficient algorithm for sampling from bridges of ergodic diffusions was proposed in Bladt and Sørensen (2007). Either of these simulation methods can be used to define the approximate likelihood function mentioned above. The main theoretical contribution of the present paper is Theorem 1, which states that, under some additional regularity conditions on the drift $\mu(\cdot, \theta)$ in (1.1), the approximation for the likelihood function obtained by the simulation of the bridge of the diffusion in (1.1) is justified because it approximates uniformly the true likelihood. For the precise statement of the result, see Section 2.

In this paper, we also propose a new algorithm for the inference of parameters when diffusion (1.1) takes a simpler form, as given in (3.1). Our method circumvents the use of numerical optimization to determine the parameters for diffusion models of the form (3.1). The estimation algorithm transforms the original problem into a related inference problem that has a unique global solution $\theta^* \in \mathbb{R}^{N+1}$ which is obtained in a single step by solving a linear system of dimension $(N+1) \times (N+1)$ given in (3.5). By Theorem 1, the related inference problem approximates uniformly on compact subsets of the parameter space the original inference problem. We also show that the approximations of the expectations that feature in linear system (3.5) converge uniformly on bounded subsets of the parameter space as the time interval between consecutive data points goes to zero (i.e., the number of observations $K+1$ goes to infinity). This property is implied by Theorem 3.

Diffusions that are not of the form (1.1) can often be transformed into the required structure by a well-known change-of-variable method [see (3.6) at the end of Section 3]. The constant diffusion coefficient requirement in (1.1) is therefore not as restrictive as it may seem at first glance. Many of the widely used diffusion processes with state-dependent volatility fall into this class. The square root process and the flexible diffusion used in Aït-Sahalia (1996) and Jones (2003b) [see (4.3) and (4.4) for the SDEs describing the models] can be dealt with in this fashion. The likelihood functions of both processes, conditional upon S&P 100 implied volatility index data, are analyzed using the EML algorithm in Section 4.2. In the case of the square root process, direct maximum likelihood estimation is also performed and it is shown that the parameter values obtained agree with the ones found using an algorithm based on the EML procedure. In this paper, we consider only one-dimensional diffusion processes, even though the EML algorithm can be applied to the multidimensional case without introducing additional computational complexity when the underlying process is reducible [see Aït-Sahalia (2008) for the precise definition]. However, extending the result of Theorem 1 to higher dimensions is a much harder problem.

The paper is organized as follows. Section 2 describes how to approximate the likelihood function and states our main theoretical result (Theorem 1). Section 3 defines and derives the EML algorithm for diffusions given by (3.1). Section 4 consists of two subsections: in Section 4.1, a comparison of the EML algorithm with exact maximum likelihood estimation and the analytic likelihood approximation method from Aït-Sahalia (2002) is performed; Section 4.2 estimates the square root and flexible diffusion processes conditional on implied volatility data. Section 5 concludes the paper. The Appendices A and B contain the proofs of Theorems 1 and 3.

## 2. The main result.

Let $x_0, \ldots, x_K$ denote $K+1$ realizations of the diffusion $X_t$ given in (1.1), observed at times $t_0, \ldots, t_K$. To avoid notational complexity, we assume evenly spaced time intervals between consecutive data points: $\Delta = t_k - t_{k-1}$ for all $k \in \{1, \ldots, K\}$. Since we are assuming that the drift $\mu(\cdot, \theta) \colon \mathbb{R} \to \mathbb{R}$ is locally Lipschitz and has linear growth, SDE (1.1) has a weakly unique solution for any starting point $x_0$ in its domain and any parameter value $\theta \in \Theta \subseteq \mathbb{R}^{N+1}$.

The starting point of our approach is the EM algorithm, which we now briefly review [see Dempster, Laird and Rubin (1977) or McLachlan and Krishnan (1997) for the general theory]. We begin by considering two consecutive data points and then apply our analysis to the entire data set. Between two consecutive observations $x_{k-1}$ and $x_k$ at times $t_{k-1}$ and $t_k$, respectively, we introduce $M-1$ evenly spaced auxiliary latent state variables $u_1, \ldots, u_{M-1}$ and define $u_0 := x_{k-1}$, $u_M := x_k$. Note that the length of the time interval between $u_{m-1}$ and $u_m$, for any $m \in \{1, \ldots, M\}$, equals $\delta := \Delta / M$. Given two observed realizations, $u_0$ and $u_M$, the task is to find the parameter $\theta = (a_0, \ldots, a_N)$ such that the value $\pi(u_M \mid u_0, \theta)$ of the conditional transition density of the diffusion $X_t$ is maximized. Consider the joint likelihood $\pi(u_M, \ldots, u_1 \mid u_0, \theta)$ of the variable $u_M$ and the latent auxiliary variables $u_1, \ldots, u_{M-1}$ conditional on $u_0$. The Markov property of the diffusion $X_t$ implies the following representation:

$$(2.1) \qquad \pi(u_M, \ldots, u_1 \mid u_0, \theta) = \prod_{m=1}^{M} \pi(u_m \mid u_{m-1}, \theta).$$

In order to formulate the EM algorithm in our setting, we need to introduce the following notation. Let random variables $U_m := X_{t_{k-1}+m\delta}$ for all $m \in \{1, \ldots, M-1\}$ correspond to the auxiliary states between the consecutive observations and let the random vector $U := (U_1, \ldots, U_{M-1})$ be the auxiliary state vector. The joint distribution of $U$ is given by the law of the bridge of the diffusion $X_t$, which starts at time $t_{k-1}$ at the level $u_0$ and finishes at the level $u_M$ at time $t_k$, denoted by $\mathbb{Q}_\theta$ (or, more accurately, $\mathbb{Q}_\theta^{\Delta, u_0, u_M}$). The subscript $\theta$ signifies the dependence of this probability law on the parameters in the model. The EM algorithm starts with some feasible value $\theta_0$ of the parameter $\theta$ and repeats the following two steps:

E-*step*: determine the conditional expectation $\theta \mapsto \mathbb{E}_{\mathbb{Q}_{\theta_n}}[\log \pi(U, u_M \mid u_0, \theta)]$;
M-*step*: maximize this expression with respect to $\theta$.

The function in the E-step of the algorithm is known as the *complete likelihood*. The important observation is that the expectation defining the complete likelihood is taken with respect to the distribution of the vector $U$, which is given by the law $\mathbb{Q}_{\theta_n}$ of the diffusion bridge. With each iteration, the value $\pi(u_M \mid u_0, \theta)$ is increased and therefore the algorithm is guaranteed to converge to a stationary point which, in some pathological cases, is not a local maximum [see McLachlan and Krishnan (1997) for convergence properties of the EM algorithm]. It is thus key to understanding the behavior of the complete likelihood $\theta \mapsto \mathbb{E}_{\mathbb{Q}_{\theta_n}}[\log \pi(U, u_M \mid u_0, \theta)]$ for any fixed parameter value $\theta_n$.

There are two problems associated with the E-step of the algorithm in our setting. The first is that the joint density $\pi(u_M, \ldots, u_1 \mid u_0, \theta)$ for the law of the process $X_t$ given by SDE (1.1) cannot be obtained in closed form. The second problem is that the law of the bridge of the diffusion $X_t$ (i.e., the process $X_t$ conditional upon $X_{t_{k-1}} = u_0$ and $X_{t_k} = u_M$) which arises in the expectation is also unknown.

Using an Euler scheme approximation for the solution of SDE (1.1), together with Markov property (2.1), one can obtain an approximation for the joint likelihood function $\pi(u_M, \ldots, u_1 \mid u_0, \theta)$ in the following way. Over any short time period of length $\delta$, the Euler scheme approximates the solution $X_{t+\delta}$ of SDE (1.1) at time $t + \delta$, conditional upon the level $X_t$, by the normal random variable $X_t + \mu(X_t, \theta)\delta + W_\delta$ with mean $X_t + \mu(X_t, \theta)\delta$ and variance $\delta$. Over a longer time period $\Delta$, a succession of such normal variables is used to approximate the original process [see Kloeden and Platen (1999), Section 10.2, for the general theory and convergence properties of Euler schemes for SDEs]. Each transition density $\pi(u_m \mid u_{m-1}, \theta)$, $m \in \{1, \dots, M\}$, in (2.1) can therefore be approximated by the normal density $\phi(u_m; u_{m-1} + \mu(u_{m-1}, \theta)\delta, \delta)$ with mean $u_{m-1} + \mu(u_{m-1}, \theta)\delta$ and variance $\delta = \frac{\Delta}{M}$ defined above. Identity (2.1) implies that an approximation for the joint likelihood $\pi(u_M, \dots, u_1 \mid u_0, \theta)$ is given by the product $\prod_{m=1}^{M} \phi(u_m; u_{m-1} + \mu(u_{m-1}, \theta)\delta, \delta)$. This approximation is useful because it depends explicitly [through the known drift function $\mu(\cdot, \theta)$] on the model parameter $\theta$ and has been used in Pedersen (1995) and Brandt and Santa-Clara (2002) to obtain an approximation for the transition density. The method, known as simulated maximum likelihood (SML), is based on the following convergence result which holds under global Lipschitz and linear growth conditions [see Theorem 2 in Pedersen (1995)]:

$$(2.2) \quad \begin{aligned} &\pi(u_M \mid u_0, \theta) \\ &= \lim_{M \to \infty} \int_{\mathbb{R}^{M-1}} \prod_{m=1}^{M} \phi(u_m; u_{m-1} + \mu(u_{m-1}, \theta)\delta, \delta) \, du_1 \cdots du_{M-1}. \end{aligned}$$

The SML algorithm uses this relationship to obtain an approximation for the likelihood function directly. As we will now show, it is possible to circumvent the difficult issue of the computation of high-dimensional integrals and obtain optimal parameter values *without* having to compute approximations of transition densities.

In the same spirit as in Jones (1998), Eraker (2001), Elerian, Chib and Shephard (2001) and Roberts and Stramer (2001), the auxiliary data points introduced at the beginning of this section are there to exploit the convergence of the discrete-time approximation to the diffusion $X_t$. As shown above, the main problem is to find the maximum of the complete likelihood $\theta \mapsto \mathbb{E}_{\mathbb{Q}_{\theta_0}}[\log \pi(U, u_M \mid u_0, \theta)]$ for any value $\theta_0$ in the parameter space. Instead of doing this, we solve an "approximate" problem, where the function $\pi(u_M, \dots, u_1 \mid u_0, \theta)$ under the expectation is replaced by the density $\prod_{m=1}^{M} \phi(u_m; u_{m-1} + \mu(u_{m-1}, \theta)\delta, \delta)$ of the Euler scheme approximation. The approximate likelihood function can then be obtained as soon as we can simulate the trajectories of the diffusion bridge [using, e.g., the algorithm in Beskos et al. (2006) or Bladt and Sørensen (2007)] governed by the law $\mathbb{Q}_{\theta_0}$. In Section 3, we will show that the approximate problem has a unique maximum that can be obtained as a solution of a linear system of size $(N+1) \times (N+1)$, where $N+1$ is the dimension of the parameter $\theta$, if the underlying diffusion takes the form (3.1).

A natural question that arises at this point concerns the quality of the approximation of the complete likelihood by the sequence of functions $(\theta \mapsto$ $\mathbb{E}_{\mathbb{Q}_{\theta_0}}[\sum_{m=1}^{M} \log \phi(U_m; U_{m-1} + \frac{\Delta}{M}\mu(U_{m-1}, \theta), \frac{\Delta}{M})])_{M \in \mathbb{N}}$. Numerical experiments in Section 4 suggest that this approximation works well. Under some additional regularity hypothesis on the drift $\mu(\cdot, \theta)$, this intuitive claim can be justified by the following theorem.

THEOREM 1. *Suppose that, in addition to the local Lipschitz condition with linear growth, we also assume that the function* $\mu: \mathbb{R} \times \mathbb{R}^{N+1} \to \mathbb{R}$ *in (1.1) is twice differentiable in the state variable with bounded second derivative. Let* $\theta_0$ *be a fixed value in the parameter space. The following equality then holds:*

$$\lim_{M \to \infty} \mathbb{E}_{\mathbb{Q}_{\theta_0}}\left[\sum_{m=1}^{M}\left(\log \pi(U_m \mid U_{m-1}, \theta)\right. \right.$$
$$\left.\left. - \log \phi\left(U_m; U_{m-1} + \frac{\Delta}{M}\mu(U_{m-1}, \theta), \frac{\Delta}{M}\right)\right)\right] = 0$$

*for all* $\theta$ *in the parameter space, where* $\phi(y; x, \delta)$ *is the normal density function with mean* $x$ *and variance* $\delta$. *Furthermore, the limit is uniform in* $\theta$ *on each compact subset of the parameter space.*

Note that the nature of Theorem 1 is fundamentally different from that of the result (2.2) above, proved in Pedersen (1995), because the expectation in the theorem is taken with respect to the law of the bridge of the diffusion $X_t$, rather than the law of the diffusion $X_t$ itself (i.e., conditional upon $X_{t_k} = u_M = x_{t_k}$ and $X_{t_{k-1}} = u_0 = x_{t_{k-1}}$). The proof of Theorem 1 is found in Appendix A. It is based on the fact that in the one-dimensional case, there exists an explicit formula for the transition density of the diffusion in terms of the Brownian bridge [see Rogers (1985)]. Because this is a special property of the one-dimensional case, the proof does not easily generalize to the multidimensional setting.

The main contribution of Theorem 1 is that it provides the theoretical basis for using the approximate complete likelihood described above with any method capable of simulating the diffusion bridge of the process defined by (1.1), including the algorithms in Beskos et al. (2006) and Bladt and Sørensen (2007), provided the regularity conditions on the drift are met. In Section 3, we use Theorem 1 to justify a key step in a new estimation algorithm for discretely observed diffusions.

**3. Expected maximum likelihood (EML) algorithm.** In this section, we are concerned with the estimation of the parameters in the SDE

$$(3.1) \quad dX_t = \mu(X_t, \theta) \, dt + dW_t \qquad \text{where } \mu(x, \theta) := g(x) + \sum_{i=0}^{N} a_i f_i(x),$$

driven by the standard Brownian motion $W_t$. The drift $\mu(\cdot, \theta) \colon \mathbb{R} \to \mathbb{R}$ is given by an arbitrary family of independent, possibly nonlinear, Lipschitz functions $g, f_i \colon \mathbb{R} \to \mathbb{R}$ with linear growth. The task is to infer the vector of coefficients $\theta := (a_0, \ldots, a_N) \in \Theta \subseteq \mathbb{R}^{N+1}$ in the drift $\mu(\cdot, \theta)$ from $K + 1$ observed realizations, $x_0, \ldots, x_K$, of the diffusion $X_t$. As we shall see, the EML algorithm consists of solving a linear system of size $(N+1) \times (N+1)$ given in (3.5) and converges to the global maximum in a single step.

Having constructed the approximation to the complete likelihood in Section 2, we now turn to the initial estimation problem. By the M-step of the EM algorithm, our task is to maximize the approximate complete likelihood function $\theta \mapsto \mathbb{E}_{\mathbb{Q}_{\theta_0}}[\sum_{m=1}^{M} \log \phi(U_m; U_{m-1} + \delta\mu(U_{m-1}, \theta), \delta)]$ for any fixed value of the model parameter $\theta_0$. The following obvious proposition is crucial to all that follows.

PROPOSITION 2. *The complete likelihood* $\theta \mapsto \mathbb{E}_{\mathbb{Q}_{\theta_0}}[\sum_{m=1}^{M} \log \phi(U_m; U_{m-1} + \delta\mu(U_{m-1}, \theta), \delta)]$ *is a nondegenerate quadratic form with a unique global maximum.*

It is clear that the complete likelihood in Proposition 2 is a nondegenerate quadratic form in $\theta$, bounded above by a constant, which implies that all of its eigenvalues must be negative. Therefore, there exists a unique global maximum.

The following simple calculation will yield the globally optimal parameter value $\theta^\star = (a_0^\star, \ldots, a_N^\star)$, which exists by Proposition 2. By setting the partial derivative with respect to each coordinate $a_j$, $j \in \{0, \ldots, N\}$, of $\theta$ in the function given in Proposition 2 to zero, we obtain the linear system $A\theta^\top = b$, where $\theta = (a_0, \ldots, a_N)$, $A = \delta \sum_{m=1}^{M} A_m$, $b = \sum_{m=1}^{M} b_m$ and, for any $m \in \{1, \ldots, M\}$,

$$(3.2) \quad A_m := \begin{pmatrix} \mathbb{E}_{\mathbb{Q}_{\theta_0}}[f_0(U_{m-1})f_0(U_{m-1})] & \cdots & \mathbb{E}_{\mathbb{Q}_{\theta_0}}[f_N(U_{m-1})f_0(U_{m-1})] \\ \vdots & \ddots & \vdots \\ \mathbb{E}_{\mathbb{Q}_{\theta_0}}[f_0(U_{m-1})f_N(U_{m-1})] & \cdots & \mathbb{E}_{\mathbb{Q}_{\theta_0}}[f_N(U_{m-1})f_N(U_{m-1})] \end{pmatrix},$$

$$(3.3) \quad \begin{aligned} b_m^\top := (\mathbb{E}_{\mathbb{Q}_{\theta_0}}[(U_m - U_{m-1} - g(U_{m-1})\delta)f_0(U_{m-1})] \cdots \\ \times \mathbb{E}_{\mathbb{Q}_{\theta_0}}[(U_m - U_{m-1} - g(U_{m-1})\delta)f_N(U_{m-1})])^\top. \end{aligned}$$

Since there exists a unique global maximum of the approximate complete likelihood, the inverse $A^{-1}$ must also exist and the unique optimal parameter value is given by $\theta^\star = (A^{-1}b)^\top$. For $K + 1$ observations of the process $X_t$, the globally optimal parameter value $\theta^\star$ is obtained in the same way. The only difference is that matrix (3.2) and vector (3.3) are computed using $M \cdot K$, rather than $M$, auxiliary and observed realizations [see (3.5)].

The globally optimal value $\theta^\star$ of the parameter vector solves the linear system whose coefficients are yet to be determined. Computing the expectations $\mathbb{E}_{\mathbb{Q}_{\theta_0}}[\cdot]$ in closed form is impossible because it requires the unknown transition density $\pi(u_{M-1}, \ldots, u_1 \mid u_M, u_0, \theta_0)$ of the bridge of the diffusion $X_t$. The key idea that helps to circumvent this problem is to replace the law of the bridge of $X_t$ with the law of the corresponding Brownian bridge in all of the coefficients of matrix (3.2) and vector (3.3). The crucial additional benefit of this substitution is that it removes the dependence of the coefficients of the linear system on the parameter $\theta_0$, which implies that the EM procedure terminates after only one iteration. Therefore, by Proposition 2, the EML algorithm is guaranteed to converge to the globally optimal parameter value $\theta^\star$ in a single step. A recent parameter estimation algorithm for general one-dimensional diffusion models given in Beskos et al. (2006), based on a sophisticated sampling method known as *retrospective sampling*, also employs the EM approach. The EM algorithm is also used in Bladt and Sørensen (2007), where the time-reversal symmetry of ergodic diffusions is exploited to sample from the corresponding bridge. Unlike in the case of the EML algorithm, in both of those settings, an iteration of E-step and M-step is required in order to obtain the stationary value of the model parameter.

We now need to consider the quality of the weak approximation of the law of the bridge of the diffusion $X_t$ (i.e., a process $X_t$ conditioned to start at $X_0 = x$ and finish at $X_\Delta = y$, where $\Delta$ is the length of the time interval between consecutive observations in the data) by the law of the Brownian bridge (i.e., a Brownian motion $W_t$ conditioned to start at $W_0 = x$ and finish at $W_\Delta = y$). This question is of importance because it tells us how far the coefficients of the linear system given by the matrix (3.2) and the target vector (3.3) are from the ones used in the EML algorithm (3.5). It is intuitively clear that when $\Delta$ goes to zero, the Brownian bridge approximation must improve in quality. Since the law of the diffusion bridge is absolutely continuous with respect to the law of the Brownian bridge, it is possible to bound the approximation error explicitly in terms of $\Delta$ and the model parameter $\theta$.

THEOREM 3. *Assume that functions $g, f_i : \mathbb{R} \to \mathbb{R}$, $i \in \{0, \ldots, N\}$, in the drift of SDE (3.1) satisfy the linear growth condition and are twice differentiable with bounded second derivatives, and let $G : \mathbb{R}^M \to \mathbb{R}$ be a polynomially bounded measurable function for some integer $M \in \mathbb{N}$. Let $\mathbb{Q}_\theta^{\Delta, x, y}$ denote the law of the bridge, starting at $x$ and finishing at $y$, of the diffusion $X_t$ that solves SDE (3.1) with the parameter value $\theta = (a_0, \ldots, a_N) \in \mathbb{R}^{N+1}$ and let $W_t$ denote the standard Brownian motion. The measure $\mathbb{Q}_\theta^{\Delta, x, y}$ is then absolutely continuous with respect to the law of the corresponding Brownian* bridge $\mathbb{W}^{\Delta,x,y}$ and the Radon–Nikodym derivative is given by

$$\frac{d\mathbb{Q}_\theta^{\Delta,x,y}}{d\mathbb{W}^{\Delta,x,y}} = \frac{L_\theta^\Delta}{\mathbb{E}_{\mathbb{W}^{\Delta,x,y}}[L_\theta^\Delta]}$$

$$where \ L_\theta^\Delta := \exp\left(-\frac{1}{2}\int_0^\Delta (\mu(W_s,\theta)^2 + \mu'(W_s,\theta))\,ds\right).$$

(a) *The following inequality holds for all $x, y$ in the domain of $X_t$ and all times $0 < t_1 < \cdots < t_M < \Delta$:*

$$|\mathbb{E}_{\mathbb{Q}_\theta^{\Delta,x,y}}[G(X_{t_1},\ldots,X_{t_M})] - \mathbb{E}_{\mathbb{W}^{\Delta,x,y}}[G(W_{t_1},\ldots,W_{t_M})]|^2$$

$$\leq \mathbb{E}_{\mathbb{W}^{\Delta,x,y}}\left[\left(\frac{L_\theta^\Delta}{\mathbb{E}_{\mathbb{W}^{\Delta,x,y}}[L_\theta^\Delta]} - 1\right)^2\right]\|G\|_2^2,$$

*where $\|G\|_2 := \mathbb{E}_{\mathbb{W}^{\Delta,x,y}}[G(W_{t_1},\ldots,W_{t_M})^2]^{1/2}$ denotes the $L^2$-norm of the random variable $G(W_{t_1},\ldots,W_{t_M})$.*

(b) *Let $S(\theta) := \sup_{z \in \mathbb{R}}\{\mu(z,\theta)^2 + \mu'(z,\theta)\}$ and $I(\theta) := \inf_{z \in \mathbb{R}}\{\mu(z,\theta)^2 + \mu'(z, \theta)\}$ be the maximum and minimum of the integrand in $L_\theta^\Delta$, respectively. The following inequality then holds:*

$$|\mathbb{E}_{\mathbb{Q}_\theta^{\Delta,x,y}}[G(X_{t_1},\ldots,X_{t_M})] - \mathbb{E}_{\mathbb{W}^{\Delta,x,y}}[G(W_{t_1},\ldots,W_{t_M})]|$$

$$\leq \frac{1}{2}\left(\exp\left(\frac{\Delta}{2}(S(\theta) - I(\theta))\right) - 1\right)\|G\|_2.$$

The absolute continuity of the measures $\mathbb{Q}_\theta^{\Delta,x,y}$ and $\mathbb{W}^{\Delta,x,y}$ is well known and the form of the Radon–Nikodym derivative in Theorem 3 follows from Lemma 1 in Beskos et al. (2006) and expressions (A.1) and (A.2) in Appendix A. The inequality in part (a) of the theorem bounds the error arising from the approximation of the law $\mathbb{Q}_\theta^{\Delta,x,y}$ by the measure $\mathbb{W}^{\Delta,x,y}$ in terms of the variance of the Radon–Nikodym derivative and the $L^2$-norm of the integrand. Since the latter is independent of the model parameter $\theta$, this inequality provides a way of bounding the error for a general integrand $G$ in terms of the second moment of the Radon–Nikodym derivative. In practice, the second moment can always be estimated by simulation, thus yielding a model-specific bound on the error of the coefficients in linear system (3.5) used in the EML algorithm. Figure 4 contains the graphs of the densities of the Radon–Nikodym derivative for the nonlinear SDE in (4.4) used in Section 4. A cursory inspection of the scale of the domains of these densities shows how tight the bound in part (a) of Theorem 3 really is, even for relatively large time steps $\Delta$.

It is intuitively clear that the Brownian bridge approximation works well for short time intervals $\Delta$ and less well as the time step grows. This view is supported by the inequality in part (b) of Theorem 3, which is a consequence of the bound in part (a). Furthermore, (b) implies that the approximation of the law $\mathbb{Q}_\theta^{\Delta,x,y}$ by $\mathbb{W}^{\Delta,x,y}$ is a good one, even for larger time steps $\Delta$, provided that the drift $\mu(\cdot,\theta)$ does not vary much as a function of the state. This implies that the method of approximation proposed in the EML algorithm would work well in the case of the diffusion with a periodic drift used in Example 1 in Beskos et al. (2006), for time steps $\Delta$ as large as $\frac{1}{2}$. Also, note that the norm of the random variable $G(W_{t_1}, \ldots, W_{t_M})$ in the Hilbert space $L^2(\mathbb{W}^{\Delta,x,y})$ is finite for a polynomially bounded function $G$ because the law $\mathbb{W}^{\Delta,x,y}$ of the Brownian bridge is Gaussian with bounded variance. The proof of Theorem 3 can be found in Appendix B.

Having replaced the law of the diffusion bridge $\mathbb{Q}_\theta^{\Delta,x,y}$ by the law of the Brownian bridge $\mathbb{W}^{\Delta,x,y}$, which is independent of the parameter $\theta$, we are left with the task of calculating the expectations in the coefficients of matrix (3.2) and vector (3.3). A numerical integration approach would be feasible because we have an explicit formula for the normal density of the marginals of the probability measure $\mathbb{W}^{\Delta,x,y}$. However, because of the numerous two-dimensional integrals in (3.3), the problem does not lend itself well to this approach.

An alternative approach is to simulate the paths of the Brownian bridge and use Monte Carlo methods to obtain the relevant expectations. This can be done by using the *modified Brownian bridge* sampler defined in Durham and Gallant (2002) and Chib and Shephard (2002), given by the following recursive formula:

$$(3.4) \qquad u_{m+1} = u_m + \frac{u_M - u_m}{M - m} + \sqrt{\frac{M - m - 1}{M - m}} \delta^{1/2} Z_m,$$

where $\delta = \Delta/M$ is the length of the time interval between consecutive auxiliary states, $u_0 = x, u_M = y$ are the initial and final points of the Brownian bridge and $Z_m \sim N(0,1)$ are independent random variables for all $m \in \{1, \ldots, M-1\}$. It is proved in Stramer and Yan (2007b) (see Proposition 1) that the joint density of the modified Brownian bridge equals the joint law of the Brownian bridge at the discretization times, which implies that scheme (3.4) introduces no discretization bias and is therefore preferable to the Euler approximation. Since the parameter $\theta$ does not appear in the evolution equation (3.4) of the modified Brownian bridge, the EML algorithm can be described as follows.

Let $x_0, \ldots, x_K$ be the $K+1$ observations at times $t_0, \ldots, t_K$ of the diffusion $X_t$ given by SDE (3.1) and let $\Delta = t_k - t_{k-1}$ for all $k \in \{1, \ldots, K\}$. Let $M-1$ be the number of the auxiliary state variables $u_{mk}$, $m \in \{1, \ldots, M-1\}$, between the observed data points $x_{k-1}$ and $x_k$ such that $u_{0k} = x_{k-1}$ and $u_{Mk} = x_k$ for all $k \in \{1, \ldots, K\}$. Let $S$ be the number of the simulations used. The EML algorithm then consists of the following simple steps.

*Step* 1. For each $k \in \{1, \ldots, K\}$ and each $s \in \{1, \ldots, S\}$, generate a Brownian bridge path $(u_{mk}^{(s)})_{m=0,\ldots,M}$ using (3.4).

*Step* 2. Find the unique solution of the linear system $A\theta^\top = b$, where

$$
(3.5)
\begin{aligned}
A_{ij} &:= \frac{\Delta}{M} \sum_{k=1}^{K} \sum_{m=1}^{M} \sum_{s=1}^{S} f_i(u_{m-1k}^{(s)}) f_j(u_{m-1k}^{(s)}) \quad \text{and} \\
b_i &:= \sum_{k=1}^{K} \sum_{m=1}^{M} \sum_{s=1}^{S} (u_{mk}^{(s)} - u_{m-1k}^{(s)} - g(u_{m-1k}^{(s)})) f_i(u_{m-1k}^{(s)})
\end{aligned}
$$

with $i, j \in \{0, \ldots, N\}$, to obtain the globally optimal parameter value $\theta^\star = (a_0^\star, \ldots, a_N^\star)$.

An appealing feature of the EML algorithm described above is that it circumvents the iterative process that is ubiquitous in the general EM framework. The invertibility of matrix $A$ is, by Proposition 2, equivalent to the nondegeneracy of the complete likelihood function, which is implied by the linear independence of the functions $f_i$ in the drift (3.1). Note that if auxiliary state variables $u_{1k}, \ldots, u_{M-1k}$, for all $k \in \{1, \ldots, K\}$, are not introduced, then we can remove the expectation operators in (3.2) and (3.3) or, equivalently, the sums over $s$ and $m$ in step 2 of the above algorithm. In this case, the EML algorithm reduces to the classic linear regression. Under the conjugate normal prior, the estimates of the drift parameters then coincide with the posterior mean in a Bayesian analysis of the coefficients.

We conclude this section with a brief comment about diffusion models with state-dependent diffusion functions. A scalar diffusion

$$dX_t = \mu(X_t, \theta) \, dt + \sigma(X_t, \vartheta) \, dW_t$$

can always be transformed using a change-of-variable $Y = \gamma(X, \vartheta) = \int^X \frac{du}{\sigma(u,\vartheta)}$, which depends on the diffusion parameter vector $\vartheta$, into

$$
(3.6)
\begin{aligned}
dY_t &= \mu_Y(X_t, \theta, \vartheta) \, dt + dW_t \\
\text{where } \mu_Y(y, \theta, \vartheta) &= \frac{\mu(\gamma^{-1}(y, \vartheta), \theta)}{\sigma(\gamma^{-1}(y, \vartheta), \vartheta)} - \frac{1}{2} \frac{\partial \sigma}{\partial x}(\gamma^{-1}(y, \vartheta), \vartheta).
\end{aligned}
$$

Note that if the original drift $\mu(\cdot, \theta)$ is affine in the parameter $\theta$, then so is the transformed drift $\mu_Y(\cdot, \theta, \vartheta)$. Therefore, an application of the EML algorithm is feasible for any fixed value of the diffusion parameter $\vartheta$. In Section 4, we will discuss how to apply the EML algorithm to diffusions of this kind [see the models given by (4.3) and (4.4)].

**4. Applications.** There are at least two potential applications for the EML algorithm. The first is the classical parameter estimation problem for diffusion models. The advantage of the EML approach is that the resulting parameter estimates are globally optimal and the bias introduced through the Euler approximation is, by Theorem 1, arbitrarily small. The second application is based on the fact that the EML algorithm is computationally very fast. The speed of the algorithm enables one to easily explore the dependence of the likelihood function on the diffusion parameter $\vartheta$ [see (3.6)], along with the dependence of the globally optimal drift parameter $\theta$ as a function of $\vartheta$ [see examples (4.3) and (4.4)]. Both of these applications will be illustrated in the present section.

4.1. *Base cases.* To test the EML algorithm for potential biases arising in the Euler and the Brownian bridge approximations, we start by establishing two base cases. The first case is an Ornstein–Uhlenbeck diffusion [see (4.1)], where the true transition density, and, therefore, the likelihood function, is known. The second case is a diffusion with a nonlinear drift [see (4.2)], where we employ the closed-form likelihood expansion from Aït-Sahalia (2002) as a benchmark, because this method is known to produce very accurate approximations of the true transition density [see, e.g., Schneider (2006), Hurn, Jeisman and Lindsay (2007), Jensen and Poulsen (2002)].

To put the EML algorithm to the test, we simulate 1000 data sets from each of the two models [the Ornstein–Uhlenbeck process in (4.1) and the nonlinear diffusion in (4.2)] with $K + 1 = 500$ observations in each data set. As mentioned above, in the first case, we perform an exact ML estimation on each of the data sets using the exact transition density in the likelihood function and in the second case, we first apply the closed-form likelihood expansion from Aït-Sahalia (2002) to obtain an approximation to the likelihood function which is then used in quasi-ML estimation. The two models are given by

$$(4.1) \qquad dX_t = (a_0 - a_1 X_t)\,dt + dW_t, \qquad\qquad \text{model A},$$

$$(4.2) \qquad dX_t = (a_0 + a_1 X_t + a_2 X_t^2)\,dt + dW_t, \qquad \text{model B},$$

with parameter values $a_0 = 10, a_1 = 2.5$ for model A and $a_0 = 1, a_1 = -1, a_2 = -0.5$ for model B. Time $t$ is measured in years and the consecutive observations in the generated data sets are one month apart. In other words, $\Delta = \frac{1}{12}$ and we choose an auxiliary state variable for each day of the month, that is, $M - 1 = 30$. To estimate the expectations in the EML algorithm, we use two sets of simulations, one containing $S = 1000$ and one $S = 200$ simulated paths. Figure 1 displays the comparison of the EML procedure using 1000 simulations with the direct ML estimation for model A. The biases and standard deviations of the parameter estimates are shown in Table

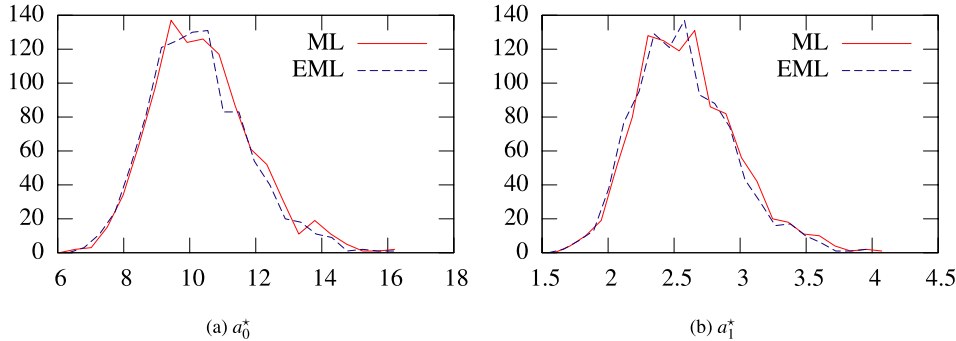

FIG. 1.  *Empirical distribution of estimates from model A. This figure shows the empirical distribution of the maximum likelihood estimator and the EML estimator over 1000 estimations for the Ornstein–Uhlenbeck model (4.1). Plotted are the maximum likelihood estimator, obtained by maximizing the true transition density using a gradient solver, and the EML estimator, obtained by using formula (3.5).*

TABLE 1

*Base cases A and B: bias and standard deviation. This table shows the estimation results for the parameters from equations (4.1) and (4.2) on 1000 data sets generated using the true transition density for model A and a very fine Euler approximation for model B (100 auxiliary data points). The benchmark estimations for model A are performed using exact maximum likelihood. The benchmark estimations for model B are obtained using closed-form likelihood expansions. The column "bias" shows the mean bias of the estimated parameters. Bias is defined as $\hat{\theta}^{(i)} - \theta_0$ for the ith estimation. The column "Std. Dev." shows the standard deviation of the parameter estimates. For the EML estimation, 30 auxiliary data points were used*

|         |       |            | **ML/Aït-Sahalia** | | **EML $S = 1000$** | | **EML $S = 200$** | |
|---------|-------|------------|--------|-----------|---------|-----------|---------|-----------|
|         |       | $\theta_0$ | Bias   | Std. Dev. | Bias    | Std. Dev. | Bias    | Std. Dev. |
| Model A | $a_0$ | 10         | 0.3902 | 1.5106    | 0.2760  | 1.4655    | 0.2765  | 1.466     |
|         | $a_1$ | 2.5        | 0.0964 | 0.3761    | 0.0678  | 0.3648    | 0.0680  | 0.3650    |
| Model B | $a_0$ | 1          | 0.1091 | 0.2769    | 0.0832  | 0.3669    | 0.0832  | 0.3669    |
|         | $a_1$ | −1         | −0.2396 | 0.5124   | −0.2352 | 0.5465    | −0.2353 | 0.5467    |
|         | $a_2$ | −0.5       | 0.0825 | 0.3525    | 0.1116  | 0.3247    | 0.1116  | 0.3246    |

1. Explicit gradients are used in the likelihood search in the case of both the exact likelihood function and the closed-form likelihood expansion.

A striking observation is the higher bias of the ML estimator over the EML estimator. The closed-form likelihood expansion and EML display similar biases. No noticeable difference can be seen by choosing 200 or 1000 simulations to approximate the expectations in the EML algorithm. In Stramer and Yan (2007a), the authors suggest that in a related problem of Monte

Carlo estimation for the transition densities of diffusions, the optimal number of simulations $S$ is of the order $M^2$, which, in the two cases discussed here, amounts to approximately 900 simulations. Also, note that the EML procedure takes about a second to produce the optimal parameter values for each of the data sets described in this subsection. The hardware used to perform these experiments was a PC with a 64-bit Xeon 2.8 MHz processor.

4.2. *Exploring the likelihood function.* The empirical features of the dynamics of equity indices such as the S&P 100 include time-varying volatility, a level effect for the volatility of the variance of return [see Jones (2003a)] and evidence for jumps [Andersen, Benzoni and Lund (2002)]. We are now going to investigate how a diffusion model, specified by a nonlinear SDE, fits the implied volatility data. In this section, we study the relation between the diffusion and drift parameters for each of the processes

$$(4.3) \quad dV_t = \kappa(\gamma - V_t)\,dt + \sigma\sqrt{V_t}\,dW_t, \qquad\qquad \text{model I,}$$

$$(4.4) \quad dV_t = (a_0 + a_1 V_t + a_2 V_t^2)\,dt + \sqrt{\sigma_1 V_t + \sigma_2 V_t^2}\,dW_t, \qquad \text{model II,}$$

conditional on S&P 100 implied volatility data given by a time series of the volatility index VXO. Empirical studies in Jones (2003a) and Aït-Sahalia and Kimmel (2007) have rejected the square root process (4.3) as a specification for the variance dynamics of S&P 100. Nevertheless, the relation between the parameters of the square root process, conditional upon real data, can be investigated[1] using the EML algorithm. The second model is a nonlinear diffusion (4.4), introduced by Aït-Sahalia (1996), and is potentially flexible enough to accommodate the rich dynamics exhibited by the S&P 100 implied volatility index data. We start by transforming the SDEs in (4.3) and (4.4) into a form with a unit diffusion coefficient using formula (3.6). For model I, we apply the transformation $y(x) = 2\sqrt{x}/\sigma$, which yields

$$(4.5) \qquad\qquad dY_t = \left(b_0 \frac{1}{Y_t} + b_1 Y_t\right)dt + dW_t.$$

Model II requires the change-of-variable $y(x) = \frac{2\log(\sqrt{x\sigma_2} + \sqrt{\sigma_1 + x\sigma_2})}{\sqrt{\sigma_2}}$, which transforms it into

$$dY_t = \mu(Y_t)\,dt + dW_t$$

---

[1] A program written in C++, which does not depend on any numerical libraries, that implements the EML algorithm in the case of the square root process, together with the VXO data used in this example, can be found at http://www.ma.ic.ac.uk/˜amijatov/Abstracts/eml.html.

with

$$\mu(y) = b_0 \frac{4}{\sqrt{e^{-2y\sqrt{\sigma_2}}(e^{2y\sqrt{\sigma_2}} - \sigma_1^2)^2/\sigma_2}} + b_1 \frac{e^{-y\sqrt{\sigma_2}}(e^{y\sqrt{\sigma_2}} - \sigma_1)^2}{\sqrt{e^{-2y\sqrt{\sigma_2}}(e^{2y\sqrt{\sigma_2}} - \sigma_1^2)^2/\sigma_2}\sigma_2}$$

$$+ b_2 \frac{e^{-2y\sqrt{\sigma_2}}(e^{y\sqrt{\sigma_2}} - \sigma_1)^4}{4\sqrt{e^{-2y\sqrt{\sigma_2}}(e^{2y\sqrt{\sigma_2}} - \sigma_1^2)^2/\sigma_2}\sigma_2^2} - \frac{e^{-y\sqrt{\sigma_2}}(\sigma_1^2 + e^{2y\sqrt{\sigma_2}})}{2\sqrt{e^{-2y\sqrt{\sigma_2}}(e^{2y\sqrt{\sigma_2}} - \sigma_1^2)^2/\sigma_2}}.$$

In the case of model I, we also perform maximum likelihood estimation with the true transition density of the square root variance (4.5).[2] The resulting parameters are $\theta = 0.0462866$, $\kappa = 5.96063$ and $\sigma = 0.455324$. Note that this parameterization ensures a stationary marginal distribution for the square root variance process.

For any fixed value of volatility $\sigma$, we can transform the time series for the volatility index using the change of coordinates $y(x) = 2\sqrt{x}/\sigma$ and perform the estimation of the parameters in (4.5) using EML. This operation takes about one second on a personal computer (64-bit Xeon 2.8 MHz, running Linux) for the given data set. By repeating this process for each value of $\sigma$ on a finite grid in the interval $[0.1, 1.2]$, we can compute the functions plotted in Figure 2(b). Using these functions, it is possible to regard the model in (4.3)

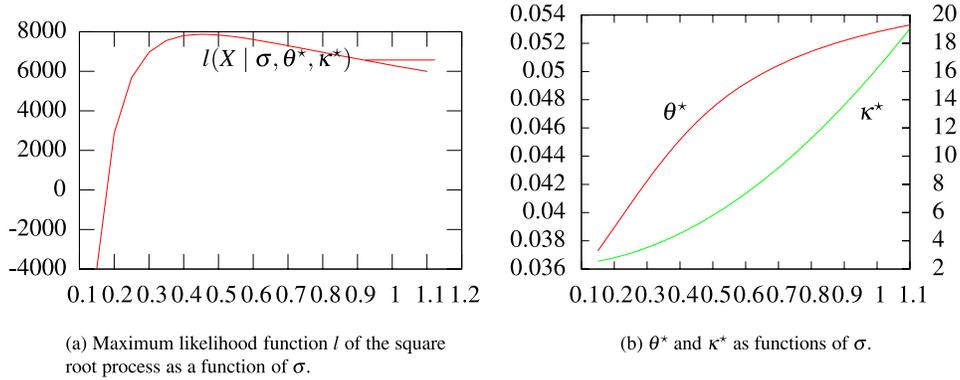

(a) Maximum likelihood function $l$ of the square root process as a function of $\sigma$.

(b) $\theta^\star$ and $\kappa^\star$ as functions of $\sigma$.

FIG. 2. *Model* I. (a) *displays the maximum likelihood function (up to a proportionality factor) of specification (4.3) as a function of $\sigma$ on the x-axis, conditional on S&P 100 implied volatility data. For a given $\sigma$, the optimal values of $\theta^\star$ and $\kappa^\star$ are computed using EML. The likelihood function is computed using the SML algorithm with the Brownian bridge importance sampler [see Durham and Gallant (2002)]. (b) displays globally optimal values $\theta^\star$ and $\kappa^\star$ as functions of $\sigma$.*

---

[2]The noncentral chi-square density is given in terms of special functions that are difficult to handle numerically for nonsymbolic computational tools with finite precision such as Fortran, C/C++ and MATLAB. We perform this estimation using Mathematica.

as having a single parameter $\sigma$ and apply the SML algorithm [as described in Durham and Gallant (2002)] or some analytic likelihood approximation to find the likelihood function for $\sigma$ [see Figure 2(a)]. This application of EML has therefore reduced the dimension of the parameter space from 3 to 1. Note that, in this approach, SML affords an additional computational advantage over other likelihood approximation methods because one can reuse the modified Brownian bridge paths that were generated for EML as an importance sampler when computing the likelihood function. Using the time series of implied volatilities, an analogous estimation can be performed for the model given by SDE (4.4).

The analysis of the likelihood functions, together with the optimal drift parameters for the processes (4.3) and (4.4), provides some interesting insights into the process specification as well as the estimation method. The first observation is that the ML estimates obtained from the true transition density of the square root process (4.3) agree with the EML estimates. This is clear from Figures 2(a) and 2(b) at the point $\sigma = 0.455324$. A visual check also reveals that even though EML is a numerical procedure, the implied globally optimal drift parameters[3] and the likelihood function are smooth in the parameters of the state-dependent volatility function. Preliminary EML estimates also suggest that the Euler approximation on a daily level is not sufficient for a nonlinear diffusion like the one in (4.4). Even a discretization of 5 subintervals per day appears too coarse. The estimates stabilize between 10 and 30 subintervals. This is in line with the findings in Roberts and Stramer (2001), Figure 4, in a similar setting. The likelihood function for model (4.4) is extremely flat in the diffusion parameters close to the optimal region (see Figure 3) and care should be taken with the estimation. Finally, a likelihood-ratio test applied to the two variance models reveals that the specification (4.4) is preferable to the square root specification.

**5. Conclusion.** This paper is concerned with an approximation procedure for the likelihood function of a nonlinear diffusion, given a discrete set of observations. The method can be used with any algorithm for sampling from the law of the diffusion bridge [e.g., Beskos et al. (2006) or Bladt and Sørensen (2007)] and is shown to converge uniformly on compact subsets of the parameter space (see Theorem 1).

We also develop a new expected maximum likelihood (EML) algorithm for the estimation of parameters governing a nonlinear diffusion process. It provides globally optimal parameter values when the drift is affine in the

---

[3]It is beneficial for the stability of the method to keep the random numbers fixed in all of the expectations arising in (3.2) and (3.3). This principle is shown to guarantee the convergence of the MCEM algorithm in Papaspiliopoulos and Sermaidis (2007).

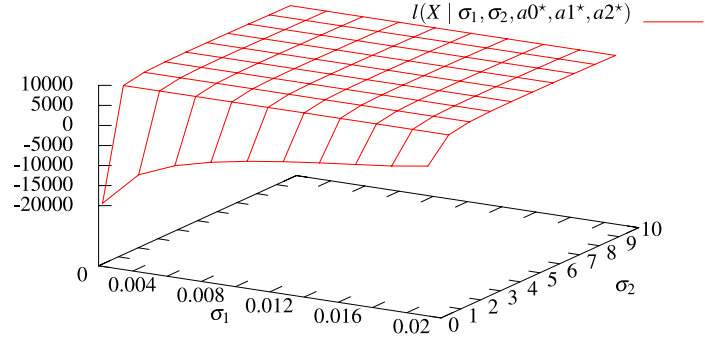

FIG. 3. *Model* II. *The likelihood of model (4.4) as a function of $\sigma_1$ (on the x-axis) and $\sigma_2$ (on the y-axis), conditional on S&P 100 implied volatility data. For given $\sigma_1$ and $\sigma_2$, the values $a_0^\star, a_1^\star$ and $a_2^\star$ are computed using EML. The likelihood function for model II given in (4.4) is computed using the SML algorithm with the Brownian bridge importance sampler [see Durham and Gallant (2002)].*

coefficients and the diffusion function is constant. For diffusions with a state-dependent volatility function, our method is used to express the likelihood as a function of the volatility parameters only, thereby significantly reduc-

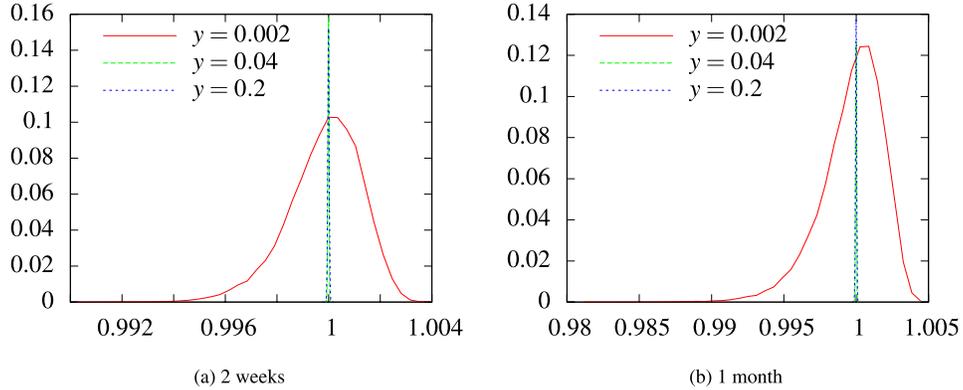

FIG. 4. *Simulated densities of the Radon–Nikodym derivative $d\mathbb{Q}_\theta^{\Delta,x,y}/d\mathbb{W}^{\Delta,x,y}$ of the law of the bridge of the diffusion $X_t$ with respect to the law of the Brownian bridge (see Theorem 3 for the precise definition). The diffusion $X_t$ used in this example is given by a nonlinear SDE (4.4). The parameter value $\theta$ is implied by the VXO data. The coordinates of $\theta$ are approximately given by $a_0 = 0.1, a_1 = -1, a_2 = -10$ and the diffusion coefficients are $\sigma_1 = 0.001$ and $\sigma_2 = 3$. The graphs correspond to time-horizons $\Delta$ equal to two weeks and one month, and a fixed starting point $x = 0.04$ for the diffusion bridge and the Brownian bridge. Three different end points $y$ of the respective bridges are chosen for each time-horizon. Note that the density in all of the cases is concentrated in a small neighborhood of 1, thus making the bound in (a) of Theorem 3 very tight, even for time intervals $\Delta$ as long as one month.*

ing the dimension of the parameter space for a gradient-based solver. The
framework is easy to implement and is guaranteed to solve the expectation
maximization problem in a single iteration. It uses auxiliary data points
and is based on two observations: the fact that the complete likelihood (i.e.,
the joint likelihood of the observed and auxiliary data points) of the Euler
scheme approximates uniformly the complete likelihood of the diffusion as
the time interval between the consecutive auxiliary data points goes to zero,
and the fact that the law of the Brownian bridge approximates well the
law of the diffusion bridge. Global optimality (Proposition 2), theoretical
bounds (Theorem 3) and asymptotic results (Theorem 1) are established,
quantifying the quality of the approximations. Additional numerical experi-
ments suggest that the method works very well for multivariate nonlinear
diffusions, even for large time intervals between observed data points.

A topic for further research is the possible extension of the EML frame-
work to the estimation of jump-diffusions. Instead of using the law of the
Brownian bridge, a semi-nonparametric density [see Gallant and Tauchen
(2009, 2006)] could be used to approximate the conditional density $p(x_{\tau_1} \mid x_{\tau_2}, x_{\tau_0}) = p(x_{\tau_2} \mid x_{\tau_1})p(x_{\tau_1} \mid x_{\tau_0})/p(x_{\tau_2} \mid x_{\tau_0}), 0 \le \tau_0 < \tau_1 < \tau_2$, of the cor-
responding bridge process with jumps. The complete likelihood function
would then be obtained by conditioning on a high-frequency discretization
(i.e., using many auxiliary state variables) of the jump-diffusion which would
identify well the jump parameters.

## APPENDIX A: PROOF OF THEOREM 1

Recall that we have a diffusion $X_t$ which is a solution of the SDE $dX_t = \mu(X_t, \theta)\, dt + dW_t$, where $\mu(\cdot, \theta) : \mathbb{R} \to \mathbb{R}$ is a bounded Lipschitz function with
bounded first and second derivatives. Without loss of generality, we can
assume that the parameter space $\Theta$ is a compact subset of $\mathbb{R}^{N+1}$. For $s > t$,
let $\pi_{t,s}(x \mid x_0, \theta)$ denote the probability density function of $X_s$ conditional
on $X_t = x_0$. It is well known that such a density exists [see (7) and (8) in
Rogers (1985)] and, by Girsanov's theorem, can be expressed as

$$\text{(A.1)} \quad \pi_{t,t+\delta}(x \mid x_0, \theta) = \frac{1}{\sqrt{2\pi\delta}} e^{-(x-x_0)^2/(2\delta)} e^{\int_{x_0}^{x} \mu(u, \theta)\, du} \Phi_\theta(\delta, x_0, x),$$

where

$$\text{(A.2)} \quad \Phi_\theta(\delta, x_0, x) := \mathbb{E}\left[\exp\left(-\frac{\delta}{2} \int_0^1 g_\theta(x_0 + u(x - x_0) + \sqrt{\delta} W_u^0)\, du\right)\right].$$

Here, $W^0$ denotes the Brownian bridge $W_u^0 := W_u - uW_1$, for $u \in [0, 1]$, and
the function $g_\theta$ is given by $g_\theta(u) := \mu'(u, \theta) + \mu(u, \theta)^2$. Our task is to prove
that the sum

$$\text{(A.3)} \quad \sum_{k=1}^{\Delta/\delta-1} \mathbb{E}_{\mathbb{Q}_{\theta_0}}\left[\log\left(\frac{\pi_{t,t+\delta}(X_{\delta(k+1)} \mid X_{\delta k}, \theta)}{\phi(X_{\delta(k+1)}; X_{\delta k} + \delta\mu(X_{\delta k}, \theta), \delta)}\right)\right]$$

converges to zero uniformly for all $\theta$ in the compact parameter space $\Theta$, where $\mathbb{Q}_{\theta_0}$ is the law of the bridge of the diffusion $X_t$, for $t \in [0, \Delta]$, which starts at $X_0 = x_0$ and finishes at $X_\Delta = x_\Delta$ for any fixed pair of real numbers $x_0, x_\Delta \in \mathbb{R}$. The function $\phi(y; x, \delta)$ in this expression is the normal density function with mean $x$ and variance $\delta$.

The first step in the proof of Theorem 1 is to understand the integrand

$$
\text{(A.4)} \quad
\begin{aligned}
\log\left( \frac{\pi_{t,t+\delta}(y \mid x, \theta)}{\phi(y; x + \delta\mu(x,\theta), \delta)} \right) = {}& \frac{\delta}{2}\mu(x,\theta)^2 - (y-x)\mu(x,\theta) \\
& + \int_x^y \mu(u, \theta)\, du + \log(\Phi_\theta(\delta, x, y)).
\end{aligned}
$$

We start with the following claim.

CLAIM 1. *The integral of the drift can be expressed as* $\int_x^y \mu(u, \theta)\, du = (y-x)\mu(x,\theta) + \frac{1}{2}(y-x)^2\mu'(x,\theta) + (y-x)^3 h_1^\theta(x,y)$, *where* $h_1^\theta(x,y)$ *is a bounded measurable function which is linear in* $\theta$.

Let $\mu(u,\theta) = \mu(x,\theta) + (u-x)\mu'(x,\theta) + \frac{1}{2}(u-x)^2\mu''(\xi_{x,u}, \theta)$ be the Taylor approximation of order 2 of the drift $\mu(\cdot, \theta)$. The point $\xi_{x,u}$ lies in the interval $(x, u) \subset \mathbb{R}$ [or in $(u, x)$, if $x$ is larger than $u$]. For any fixed $y$ in $\mathbb{R}$, we can integrate this representation of $\mu(\cdot, \theta)$ to obtain the representation of the integral which is given in Claim 1. We need to check that the function $h_1^\theta(x,y) := \frac{1}{(y-x)^3}\int_x^y (u-x)^2 b^\theta(u,x)\, du$, for $x \neq y$, is bounded and measurable. Here, $b^\theta(u,x)$ is given by the quotient $b^\theta(u,x) := 2\frac{\mu(u,\theta) - \mu(x,\theta) - (u-x)\mu'(x,\theta)}{(u-x)^2}$ if $u \neq x$ and is zero otherwise. Note that the function $b^\theta$ is bounded since $\mu''(\cdot, \theta)$ is bounded and linear in $\theta$. Since the set $\Theta$ is compact, the estimate $|h_1^\theta(x,y)| \leq \frac{C}{(y-x)^3}\int_x^y (u-x)^2\, du = \frac{C}{3}$ holds for all $y > x$ and some constant $C$. A similar bound holds for $y < x$. It follows from the definition of $b^\theta$ that it is measurable on $\mathbb{R} \times \mathbb{R}$ since it is continuous outside the zero measure set $\{(x,x) : x \in \mathbb{R}\}$. Fubini's theorem implies that $h_1^\theta(x,y)$ must therefore also be measurable. This proves Claim 1.

The next task is to relate the asymptotic behavior of the function $\log(\Phi_\theta(\delta, x, y))$ to the drift $\mu(\cdot, \theta)$. This will be achieved in Claims 2 and 3.

CLAIM 2. *There exist constants* $\delta_0 > 0$ *and* $A_0 > 0$ *such that for all* $\delta \in [0, \delta_0]$, *we have the following equality:*

$$
\log(\Phi_\theta(\delta, x, y)) + \frac{\delta}{2}\int_0^1 g_\theta(x + u(y-x))\, du = \delta^{3/2}A^\theta(x, y),
$$

*where* $A^\theta(x,y)$ *is a measurable function of* $x$ *and* $y$ *that satisfies* $|A^\theta(x,y)| < A_0$ *for all* $x, y \in \mathbb{R}$ *and all* $\theta \in \Theta$.

By Lagrange's mean value theorem, we obtain the following equalities:

$$g_\theta(x + u(y-x) + \sqrt{\delta}W_u^0) - g_\theta(x + u(y-x))$$
$$= \sqrt{\delta}W_u^0 g_\theta'(x + u(y-x) + \xi)$$
$$= \sqrt{\delta}W_u^0 X^\theta(x, y, u),$$

where the real number $\xi \in (0, \sqrt{\delta}W_u^0)$ depends on $x, y, \delta$ and $W_u^0$. Note that the random variable $X^\theta(x, y, u)$ is measurable since it coincides with the quotient $[g_\theta(x + u(y-x) + \sqrt{\delta}W_u^0) - g_\theta(x + u(y-x))]/(\sqrt{\delta}W_u^0)$ on a complement of the set $\{W_u^0 = 0\}$, which has probability zero. This representation also implies that for any fixed path of $W_u^0$, the function $X^\theta(x, y, u)$ is jointly measurable in the variables $x, y$ and $u$ and that it is quadratic in the parameter $\theta$. The assumptions on the drift $\mu(\cdot, \theta)$ and the above equality imply that the random variable $X^\theta(x, y, u)$ is bounded for all possible triplets $(x, y, u)$.

Let the variable $Y^\theta(x, y)$ denote the integral $Y^\theta(x, y) := \int_0^1 X^\theta(x, y, u) W_u^0 \, du$ and let $Z^\theta(x, y)$ denote the random variable $Z^\theta(x, y) := \frac{1}{2} \sum_{k=1}^\infty (-\frac{\delta^{3/2}}{2})^{k-1} \times \frac{Y^\theta(x,y)^k}{k!}$. The left-hand side of the equality in Claim 2 can now be rewritten as

$$(A.5) \qquad \log \mathbb{E}\left[ \exp\left( -\frac{\delta^{3/2}}{2} Y^\theta(x, y) \right) \right] = \log \mathbb{E}[1 - \delta^{3/2} Z^\theta(x, y)].$$

From the definition of $W_u^0$ and the fact that $X^\theta(x, y, u)$ is bounded by a positive constant $C$ independent of $\theta$, we can see that $|Y^\theta(x, y)| \leq C(\int_0^1 |W_u| \, du + |W_1|)$. This implies the bound $|Z^\theta(x, y)| \leq e^{|Y^\theta(x,y)|} \leq e^{C \int_0^1 |W_u| \, du} e^{C|W_1|}$, which, combined with the Cauchy–Schwarz inequality, yields

$$(A.6) \qquad \mathbb{E}[|Z^\theta(x, y)|] \leq \mathbb{E}[e^{2C \int_0^1 |W_u| \, du}]^{1/2} \mathbb{E}[e^{2C|W_1|}]^{1/2}.$$

This means that the expectation of $Z^\theta(x, y)$ exists and is bounded above and below uniformly in $x$ and $y$ for all parameters $\theta \in \Theta$. Therefore, there exists $\delta_0 > 0$ such that $-\frac{1}{2} < \delta \mathbb{E}[Z^\theta(x, y)] < \frac{1}{2}$ holds for all $x, y$, if $\delta \in [0, \delta_0]$.

The right-hand side of (A.5) can now be expressed as

$$\log \mathbb{E}[1 - \delta^{3/2} Z^\theta(x, y)] = \log(1 - \delta^{3/2} \mathbb{E}[Z^\theta(x, y)])$$
$$= \delta^{3/2} \sum_{k=1}^\infty (\delta^{3/2})^{k-1} \frac{\mathbb{E}[Z^\theta(x, y)]^k}{k}$$

since the series $\log(1 - z) = \sum_{k=1}^\infty z^k/k$ converges uniformly on compact subsets of $(-1, 1)$. We can define $A^\theta(x, y) := \sum_{k=1}^\infty (\delta^{3/2})^{k-1} \mathbb{E}[Z^\theta(x, y)]^k/k$, which is therefore measurable and uniformly bounded for all $x, y$ and $\delta \in [0, \delta_0]$. It also follows from the bound in (A.6) that there exists a constant $A_0 > 0$ such that $|A^\theta(x, y)| < A_0$ for all $x, y$ and all parameter values $\theta \in \Theta$. This proves Claim 2.

CLAIM 3. *The equality $\int_0^1 g_\theta(x + u(y-x)) \, du = g_\theta(x) + (y-x)h_2^\theta(x,y)$ holds for all $x$ and $y$ and some bounded measurable function $h_2^\theta(x,y)$ which is quadratic in the parameter $\theta$.*

By Lagrange's mean value theorem, we obtain $g_\theta(x + u(y-x)) = g_\theta(x) + g_\theta'(\xi)u(y-x)$, where $|\xi - x| < u|y - x|$. It is clear that $[g_\theta(x + u(y-x)) - g_\theta(x)]/(y - x)$ is a measurable function defined on $((\mathbb{R} \times \mathbb{R}) - \{(x,x) : x \in \mathbb{R}\}) \times [0,1]$. Since the diagonal $\{(x,x) : x \in \mathbb{R}\}$ has Lebesgue measure zero in $\mathbb{R} \times \mathbb{R}$, we can extend the quotient to a bounded measurable function on $\mathbb{R} \times \mathbb{R} \times [0,1]$. If we integrate the above equality over the interval $[0,1]$, we obtain $\int_0^1 g_\theta(x + u(y-x)) \, du = g_\theta(x) + (y-x)\int_0^1 g_\theta'(\xi)u \, du$. It follows from Fubini's theorem that the last integral, denoted by $h_2^\theta(x,y)$, is a measurable function of $x$ and $y$. It is also clear that $h_2^\theta(x,y)$ is bounded, since the integrand $(x,y,u) \mapsto g_\theta'(\xi)u$ is bounded on its domain, and that $h_2^\theta(x,y)$ is quadratic in the parameter $\theta$. This completes the proof of Claim 3.

We can now apply Claims 1, 2 and 3 to the equality in (A.4) to obtain the following representation:

$$\log\left(\frac{\pi_{t,t+\delta}(y \mid x, \theta)}{\phi(y; x + \delta\mu(x,\theta), \delta)}\right)$$
$$= \frac{1}{2}\mu'(x,\theta)[(y-x)^2 - \delta] + (y-x)^3 h_1^\theta(x,y)$$
$$- \frac{1}{2}\delta h_2^\theta(x,y)(y-x) + \delta^{3/2}A^\theta(x,y),$$

where $h_1^\theta(x,y)$, $h_2^\theta(x,y)$ and $A^\theta(x,y)$ are bounded measurable functions. The sum in (A.3) can now be decomposed naturally into four sums. The first three will tend to zero as $\delta$ goes to zero, by Lemma 4. The fourth one can easily be bounded as follows:

$$\delta^{3/2}\left|\sum_{k=1}^{\Delta/\delta-1} \mathbb{E}_{\mathbb{Q}_{\theta_0}}[A^\theta(X_{\delta k}, X_{\delta(k+1)})]\right| \leq \delta^{3/2}\left(\frac{\Delta}{\delta} - 1\right)A_0$$
$$= \sqrt{\delta}(\Delta - \delta)A_0,$$

where the constant $A_0$ is as in Claim 2 and hence converges to zero uniformly in $\theta$. Also, note that since the functions $\mu'(x,\theta)$, $h_1^\theta(x,y)$ and $h_2^\theta(x,y)$ are at most quadratic in the parameter $\theta$ which takes values in a compact region $\Theta$, it is enough to state and prove Lemma 4 for functions that do not depend on $\theta$ and still obtain uniform convergence in $\theta$. The proof of Theorem 1 will therefore be complete as soon as we prove the following lemma.

LEMMA 4. *Let $f_\delta : \mathbb{R} \times \mathbb{R} \to \mathbb{R}$ be one of the following functions:*

(1) $f_\delta(x, y) = h(x)((y - x)^2 - \delta)$ *for any bounded differentiable function* $h$ *with bounded first derivative;*

(2) $f_\delta(x, y) = h(x, y)(y - x)^3$, *where* $h(x, y)$ *is a bounded measurable function; or*

(3) $f_\delta(x, y) = \delta h(x, y)(y - x)$, *where* $h(x, y)$ *is as in (2).*

*The following then holds:*

$$\lim_{\delta \to 0} \sum_{k=1}^{\Delta/\delta - 1} \mathbb{E}_{\mathbb{Q}_{\theta_0}}[f_\delta(X_{k\delta}, X_{(k+1)\delta})] = 0.$$

PROOF. Let $A(\delta)$ denote the sum in the above limit. Since the diffusion $X_t$ is a Markov process, we can express $A(\delta)$ as

$$(A.7) \qquad A(\delta) = \sum_{k=1}^{\Delta/\delta - 1} \mathbb{E}_{\mathbb{Q}_{\theta_0}}[\mathbb{E}_{\mathbb{Q}_{\theta_0}}[f_\delta(X_{k\delta}, X_{(k+1)\delta})|X_{k\delta}, X_\Delta]].$$

Conditional densities, which arise in the expectations in (A.7), can be expressed using the formula in (A.1) for the transition density of $X_t$. We can therefore rewrite $A(\delta)$ as

$$A(\delta) = \sum_{k=1}^{\Delta/\delta - 1} \int_\mathbb{R} \frac{\pi_{0,k\delta}(x|x_0)\pi_{k\delta,\Delta}(x_\Delta|x)}{\pi_{0,\Delta}(x_\Delta|x_0)} dx$$
$$\times \int_\mathbb{R} f_\delta(x, y) \frac{\pi_{k\delta,(k+1)\delta}(y|x)\pi_{(k+1)\delta,\Delta}(x_\Delta|y)}{\pi_{k\delta,\Delta}(x_\Delta|x)} dy.$$

To simplify the notation, let $N := \frac{\Delta}{\delta} - 1$. Note that we are always choosing $\delta$ so that $N$ is an integer. The above expression for $A(\delta)$ implies that the lemma will follow if we prove the following equalities:

$$\lim_{N \to \infty} \sum_{k=1}^{N-1} \frac{1}{\delta^{3/2}\sqrt{k(N-k)}}$$
$$\times \int_{\mathbb{R} \times \mathbb{R}} f_\delta(x, y) \exp\Big(-\frac{(x_\Delta - y)^2}{2\Delta((N-k)/(N+1))}$$
$$-\frac{(y - x)^2}{2\Delta/(N+1)} - \frac{(x - x_0)^2}{2\Delta(k/(N+1))}\Big)$$
$$\times \Phi\Big(\Delta\frac{N-k}{N+1}, y, x_\Delta\Big)\Phi\Big(\frac{\Delta}{N+1}, x, y\Big)$$
$$(A.8)$$
$$\times \Phi\Big(\Delta\frac{k}{N+1}, x_0, x\Big) dx\, dy = 0,$$

$$\lim_{\delta \to 0} \int_{\mathbb{R}} f_\delta(x, x_\Delta) \frac{1}{\sqrt{\delta(\Delta - \delta)}} \exp\left( -\frac{(x_\Delta - x)^2}{2\delta} - \frac{(x - x_0)^2}{2(\Delta - \delta)} \right)$$

(A.9)

$$\times \, \Phi(\delta, x, x_\Delta) \Phi(\Delta - \delta, x_0, x) \, dx = 0.$$

Condition (A.9) corresponds to the last summand in (A.7), while condition (A.8) accounts for the rest of the sum in (A.7).

Let us start by proving (A.9). Note that since, by assumption, the drift $\mu$ is bounded and has a bounded first derivative, the function $\Phi$ must be bounded (as mentioned above, we are omitting dependence on the parameter $\theta_0$ as it is fixed throughout). The function $h$ in the definition of $f_\delta$ is also bounded, so the absolute value of the integral in (A.9) is smaller than $\delta C \int_{\mathbb{R}} |v^2 - 1| e^{-v^2/2} \, dv$, $\delta^{3/2} C \int_{\mathbb{R}} |v^3| e^{-v^2/2} \, dv$ and $\delta^{3/2} C \int_{\mathbb{R}} |v| e^{-v^2/2} \, dv$ for $f_\delta$ given by (1), (2) and (3), respectively. These bounds are obtained by the change of variable $v = (x_\Delta - x)/\sqrt{\delta}$ in the integral in (A.9). In each of the three cases, the constant $C$ is independent of $\delta$. This proves (A.9).

We are now left with the harder problem of proving (A.8). Note that the integrand is absolutely integrable over $\mathbb{R} \times \mathbb{R}$, which, by Fubini's theorem, implies that we are free to choose any order of integration in the double integral. We will first prove case (1), where $f_\delta(x, y) = h(x)((y - x)^2 - \delta)$. This case is harder than (2) and (3), which will be dealt with at the end of the proof.

By substituting $u = (x - y)/\sqrt{\delta}$ and integrating over the state variable $x$ first, we can rewrite the sum in (A.8) as

(A.10) $\quad B(u, N) := \int_{\mathbb{R}} (u^2 - 1) e^{-u^2/2} C(u, N) \, du,$

(A.11) $\quad C(u, N) := \sum_{k=1}^{N-1} \frac{F(u, k, N)}{\sqrt{k(N-k)}},$

$$F(u, k, N) := \int_{\mathbb{R}} h\left( u\sqrt{\frac{\Delta}{N+1}} + y \right)$$

$$\times \exp\left( -\frac{(u\sqrt{\Delta/(N+1)} + y - x_0)^2}{2\Delta(k/(N+1))} \right.$$

(A.12)

$$\left. -\frac{(x_\Delta - y)^2}{2\Delta((N-k)/(N+1))} \right)$$

$$\times \Phi\left( \Delta \frac{N-k}{N+1}, y, x_\Delta \right) \Phi\left( \frac{\Delta}{N+1}, u\sqrt{\frac{\Delta}{N+1}} + y, y \right)$$

$$\times \Phi\left( \Delta \frac{k}{N+1}, x_0, u\sqrt{\frac{\Delta}{N+1}} + y \right) dy.$$

We shall now prove the following statements, which will imply the lemma:

(a) $\lim_{N\to\infty} \sum_{k=1}^{N-1} \frac{1}{\sqrt{k(N-k)}} = \pi$;

(b) $\lim_{N\to\infty} C(u,N)$ exists for every $u \in \mathbb{R}$;

(c) $\lim_{N\to\infty} C(u,N)$ is independent of $u$;

(d) $\lim_{N\to\infty} B(u,N) = \int_{\mathbb{R}}(u^2-1)e^{-u^2/2}(\lim_{N\to\infty} C(u,N))\,du$.

Note that (c) and (d) together imply that $\lim_{N\to\infty} B(u,N) = (\lim_{N\to\infty} C(u, N)) \int_{\mathbb{R}}(u^2-1)e^{-u^2/2}\,du = 0$, which proves the lemma in the case where $f_\delta$ is given by (1). In order to prove (a), we need to consider the graph of the mapping $x \mapsto 1/\sqrt{x(N-x)}$ in the interval $(0,N)$. We can, without loss of generality, assume that $N$ is an even integer by choosing the decreasing time interval $\delta$ correspondingly [recall that $\delta = \Delta/(N+1)$]. By inspecting the area under the graph of this map, it is easy to see that the following inequalities must hold:

$$\int_1^{N/2} \frac{dx}{\sqrt{x(N-x)}} + \frac{2}{N} + \int_{N/2}^{N-1} \frac{dx}{\sqrt{x(N-x)}}$$

$$\leq \sum_{k=1}^{N-1} \frac{1}{\sqrt{k(N-k)}}$$

$$\leq \int_0^{N/2} \frac{dx}{\sqrt{x(N-x)}} + \int_{N/2+1}^N \frac{dx}{\sqrt{x(N-x)}}.$$

It is easy to check that $\int_0^x \frac{dz}{\sqrt{z(y-z)}} = 2\arctan(\sqrt{\frac{x}{y-x}})$ for all $x \in (0,y)$, where $y$ is any positive real number. Using this formula to calculate the integrals in the above inequalities, we obtain the following relations:

$$2(\arctan(\sqrt{N-1}) - \arctan(1/\sqrt{N-1}) + 2/N)$$

$$\leq \sum_{k=1}^{N-1} \frac{1}{\sqrt{k(N-k)}}$$

$$\leq 2\left(\arctan(1) + \frac{\pi}{2} - \arctan\left(\sqrt{\frac{N/2+1}{N/2-1}}\right)\right)$$

for every even integer $N$. It is now clear that the limit of the sum in (a) equals $\pi$.

We are now going to prove (b) and (c). Note that we can assume, without loss of generality, that the function $h$ in (A.12) is nonnegative [i.e., $h(y) \geq 0$ for all $y \in \mathbb{R}$] because we can express it as $h = h_+ - h_-$, where $h_+(y) = \max\{h(y), 0\}$ and $h_-(y) = \max\{-h(y), 0\}$ are nonnegative functions. Note that the functions $h_+$ and $h_-$ are nondifferentiable on a set which is at most countable[4] and therefore has Lebesgue measure zero.

We will start by showing that the limit $\lim_{N\to\infty} C(0,N)$ exists. Since $h$ is nonnegative, the same is true of $F(0,k,N)$. Using the same reasoning as above, we obtain the following inequalities:

$$
\begin{aligned}
(A.13) \qquad & \sum_{k=1}^{N/2-1} F(0,k,N) \int_k^{k+1} \frac{dx}{\sqrt{x(N-x)}} \\
& + F\left(0, \frac{N}{2}, N\right) \frac{2}{N} + \sum_{k=N/2}^{N-2} F(0,k+1,N) \int_k^{k+1} \frac{dx}{\sqrt{x(N-x)}} \\
& \leq \sum_{k=1}^{N-1} \frac{F(0,k,N)}{\sqrt{k(N-k)}} \\
& \leq \sum_{k=0}^{N/2-1} F(0,k+1,N) \int_k^{k+1} \frac{dx}{\sqrt{x(N-x)}} \\
& + \sum_{k=N/2+1}^{N-1} F(0,k,N) \int_k^{k+1} \frac{dx}{\sqrt{x(N-x)}}.
\end{aligned}
$$

Since the function $x \mapsto 1/\sqrt{x(N-x)}$ is integrable on the interval $[0,N]$ and $F$ is bounded, the inequalities in (A.13), together with the next claim, imply that the limit $\lim_{N\to\infty} C(0,N)$ exists.

CLAIM.

$$
\begin{aligned}
\lim_{N\to\infty} & \sum_{k=1}^{N/2-1} |F(0,k+1,N) - F(0,k,N)| \int_k^{k+1} \frac{dx}{\sqrt{x(N-x)}} \\
& = \lim_{N\to\infty} \sum_{k=N/2+1}^{N-2} |F(0,k+1,N) - F(0,k,N)| \int_k^{k+1} \frac{dx}{\sqrt{x(N-x)}} = 0.
\end{aligned}
$$

---

[4] The set $A$, where $h_+$ is not differentiable, can be expressed in terms of $h$ as follows: $A = h^{-1}(0) \cap (h'^{-1}(0,\infty) \cup h'^{-1}(-\infty,0))$. This is because $h_+$ is nondifferentiable precisely at those zeros of $h$ where the derivative $h'$ is nonzero. Since $h'$ is continuous, the set $h'^{-1}(0,\infty)$ [resp., $h'^{-1}(-\infty,0)$] is open in $\mathbb{R}$ and can therefore be expressed as a disjoint union of open intervals. Recall that there can only be countably many disjoint open intervals in $\mathbb{R}$. Since function $h$ is monotonic on each of these intervals, there can be at most one zero in each interval. Therefore, the set $A$ of such zeros of $h$ must be countable.

The two cases in the claim are very similar and can be deduced using the same methodology. We will give details only for the first case. Let $k \in \{1, \ldots, \frac{N}{2} - 1\}$ and recall that $\Phi$ and $h$ are bounded functions with bounded derivatives. We can therefore obtain, using (A.12), the following estimate:

$$(A.14) \qquad |F(0, k+1, N) - F(0, k, N)| \leq \sum_{i=1}^{4} C_i I_i(k, N),$$

where $C_i$, $i = 1, \ldots, 4$, are constants independent of $k$ and $N$, and $I_i(k, N)$, $i = 1, \ldots, 4$, are integrals given by:

$$I_1(k, N) := \int_{\mathbb{R}} \left| \exp\left( -\frac{(x_\Delta - y)^2}{2\Delta((N - k - 1)/(N + 1))} \right) \right.$$
$$\left. - \exp\left( -\frac{(x_\Delta - y)^2}{2\Delta((N - k)/(N + 1))} \right) \right| dy;$$

$$I_2(k, N) := \int_{\mathbb{R}} \left| \exp\left( -\frac{(y - x_0)^2}{2\Delta((k + 1)/(N + 1))} \right) \right.$$
$$\left. - \exp\left( -\frac{(y - x_0)^2}{2\Delta(k/(N + 1))} \right) \right|$$
$$\times \exp\left( -\frac{(x_\Delta - y)^2}{2\Delta((N - k)/(N + 1))} \right) dy;$$

$$I_3(k, N) := \int_{\mathbb{R}} \left| \Phi\left( \Delta\frac{N - k - 1}{N + 1}, y, x_\Delta \right) - \Phi\left( \Delta\frac{N - k}{N + 1}, y, x_\Delta \right) \right|$$
$$\times \exp\left( -\frac{(x_\Delta - y)^2}{2\Delta((N - k)/(N + 1))} \right) dy;$$

$$I_4(k, N) := \int_{\mathbb{R}} \left| \Phi\left( \Delta\frac{k + 1}{N + 1}, x_0, y \right) - \Phi\left( \Delta\frac{k}{N + 1}, x_0, y \right) \right|$$
$$\times \exp\left( -\frac{(x_\Delta - y)^2}{2\Delta((N - k)/(N + 1))} \right) dy.$$

We can now estimate the integrals $I_i(k, N)$, $i = 1, \ldots, 4$, for any $k \in \{1, \ldots, \frac{N}{2} - 1\}$. For $I_1(k, N)$, we get

$$I_1(k, N) \leq \int_{\mathbb{R}} \exp\left( -\frac{(x_\Delta - y)^2}{\Delta(N/(N + 1))} \right)$$
$$\times \left| 1 - \exp\left( -\frac{(x_\Delta - y)^2}{2\Delta} \frac{N + 1}{(N - k - 1)(N - k)} \right) \right| dy$$
$$\leq \int_{\mathbb{R}} \exp(-(x_\Delta - y)^2/\Delta)$$

$$\times \left(1 - \exp\left(-\frac{(x_\Delta - y)^2}{\Delta}\frac{N+1}{N(N/2-1)}\right)\right) dy.$$

The last integral in this inequality converges to zero, by the dominated convergence theorem. Therefore, $I_1(k, N)$ goes to zero with increasing $N$ uniformly in $k$.

In order to bound $I_2(k, N)$, note that the function $x \mapsto \exp(-(y-x_0)^2/(2\Delta x))$ has a derivative which is bounded (for each $y \in \mathbb{R}$) on the interval $x \in [0, 1/2]$. Since $k \in \{1, \ldots, N/2 - 1\}$, the value $x = k/(N+1)$ lies in the interval $(0, 1/2)$. For $k$ such that $k/(N+1) \le (N+1)^{-1/4}$, we get the following bound:

$$
\begin{aligned}
\text{(A.15)} \quad I_2(k, N) &\le \int_\mathbb{R} \left(\exp\left(-\frac{(y-x_0)^2}{2\Delta 2(N+1)^{-1/4}}\right)\right. \\
&\qquad\qquad \left. + \exp\left(-\frac{(y-x_0)^2}{2\Delta(N+1)^{-1/4}}\right)\right) dy \\
&= \frac{\sqrt{2\pi\Delta}(1+\sqrt{2})}{(N+1)^{1/8}},
\end{aligned}
$$

which is uniform in all $k$ that satisfy the above condition.

For $k \in \{1, \ldots, N/2 - 1\}$ such that $k/(N+1) > (N+1)^{-1/4}$, the Lagrange mean value theorem implies the existence of some $\xi$ in the interval $(\frac{k}{N+1}, \frac{k+1}{N+1})$ such that

$$
\begin{aligned}
&\left|\exp\left(-\frac{(y-x_0)^2}{2\Delta(k/(N+1)+1/(N+1))}\right) - \exp\left(-\frac{(y-x_0)^2}{2\Delta(k/(N+1))}\right)\right| \\
&\qquad \le \frac{1}{N+1}\frac{(y-x_0)^2}{2\Delta\xi^2}e^{-(y-x_0)^2/(2\Delta\xi)} \\
&\qquad \le \frac{1}{N+1}\frac{(y-x_0)^2\sqrt{N+1}}{2\Delta}.
\end{aligned}
$$

The last inequality is independent of $k$ and holds for every $y \in \mathbb{R}$. It therefore implies that

$$I_2(k, N) \le \frac{C}{\sqrt{N+1}}$$

for all $k \in \{1, \ldots, \frac{N}{2} - 1\}$ which satisfy the inequality $\frac{k}{N+1} > (N+1)^{-1/4}$ (the constant $C$ in the above expression is independent of $k$). This, together with the inequality in (A.15), implies that $I_2(k, N)$ converges to zero uniformly in $k$.

It is shown in Rogers (1985) (see page 160) that the partial derivative $\frac{\partial \Phi}{\partial \delta}$ of the function $\Phi(\delta, x_0, x)$ defined in (A.2) equals the expectation of the

derivative

$$\frac{\partial}{\partial\delta}\left[\exp\left(-\frac{1}{2}\delta\int_0^1 g(x_0 + u(x-x_0) + \sqrt{\delta}W_u^0)\,du\right)\right].$$

Using the fact that $g$ is bounded and has a bounded first derivative, we conclude that $\frac{\partial\Phi}{\partial\delta}$ is also bounded. Therefore, we can apply Lagrange's mean value theorem to prove that $I_3(k, N)$ and $I_4(k, N)$ converge to zero uniformly in $k$ as $N$ goes to infinity.

We have just shown that $\lim_{N\to\infty} M(N) = 0$, where $M(N) := \max\{|F(0, k+1, N) - F(0, k, N)| : k = 1, \ldots, \frac{N}{2} - 1\}$. The inequalities

$$0 \le \sum_{k=1}^{N/2-1} |F(0, k+1, N) - F(0, k, N)| \int_k^{k+1} \frac{dx}{\sqrt{x(N-x)}}$$

$$\le M(N) \int_1^{N/2} \frac{dx}{\sqrt{x(N-x)}},$$

together with (a) and the fact that $M(N)$ converges to zero, prove the claim. We have therefore proven that the limit $\lim_{N\to\infty} C(0, N)$ exists.

In order to complete the proof of (b) and (c), we need to understand the behavior of

$$\frac{\partial C}{\partial u}(u, N) = \sum_{k=1}^{N-1} \frac{1}{\sqrt{k(N-k)}} \frac{\partial F}{\partial u}(u, k, N)$$

for all $u \in \mathbb{R}$. The key observation here is that every instance of the parameter $u$ in definition (A.12) of the function $F(u, k, N)$ is of the form $\frac{u}{\sqrt{N+1}}$. It is therefore natural to expect that the partial derivative $\frac{\partial F}{\partial u}(u, k, N)$ tends to zero with increasing $N$ for every fixed value of $u$. This is precisely what we shall now prove.

It follows from (A.12) that the inequality

$$(A.16) \qquad\qquad \left|\frac{\partial F}{\partial u}(u, k, N)\right| \le \frac{D_0}{\sqrt{N+1}} \sum_{i=1}^4 J_i$$

holds, where $D_0$ is a positive constant independent of $u$, $k$ and $N$, and the integrals $J_i$, $i = 1, \ldots, 4$, are of the following form:

$$J_1 := \int_{\mathbb{R}} \left|\Phi_2\left(\Delta\frac{k}{N+1}, x_0, u\sqrt{\frac{\Delta}{N+1}} + y\right)\right| \exp\left(-\frac{(x_\Delta - y)^2}{2\Delta}\right) dy;$$

$$J_2 := \int_{\mathbb{R}} \left|h'\left(u\sqrt{\frac{\Delta}{N+1}} + y\right)\right| \exp\left(-\frac{(x_\Delta - y)^2}{2\Delta}\right) dy;$$

$$J_3 := \int_{\mathbb{R}} \left| \Phi_1 \left( \frac{\Delta}{N+1}, u\sqrt{\frac{\Delta}{N+1}} + y, y \right) \right| \exp\left( -\frac{(x_\Delta - y)^2}{2\Delta} \right) dy;$$

$$J_4 := \int_{\mathbb{R}} \frac{|u\sqrt{\Delta/(N+1)} + y - x_0|}{\Delta(k/(N+1))} \exp\left( -\frac{(u\sqrt{\Delta/(N+1)} + y - x_0)^2}{2\Delta(k/(N+1))} \right) dy.$$

The functions $\Phi_1$ and $\Phi_2$ denote the derivatives of the function $\Phi$ given in (A.2) with respect to the first and second state variable, respectively. It is shown in Rogers (1985) that these derivatives exist and that they are bounded.

In order to obtain the bound in (A.16), we had to exchange the order of differentiation and integration [see the definition of function $F$ in (A.12)]. This can be justified by the dominated convergence theorem since the difference quotients are bounded above by the function

$$y \mapsto \sup_{u' \in (u-1, u+1)} \left[ \frac{\partial f}{\partial u}(u', y) \right] \exp\left( -\frac{(x_\Delta - y)^2}{2\Delta} \right)$$

for each $u \in \mathbb{R}$, where $f$ is the integrand in (A.12). This function is clearly in $L^1(\mathbb{R})$ and the dominated convergence theorem applies. Also, note that the integral $J_2$ is well defined because, as noted earlier (see page 23), the function $h$ is nondifferentiable only on a set of measure zero.

We can now proceed to estimate the integrals $J_i$, $i = 1, \dots, 4$. Since the functions $\Phi_1, \Phi_2$ and $h'$ are bounded, the integrals $J_1$, $J_2$ and $J_3$ are also bounded above by a constant for all $u$, $N$ and $k$. By introducing a change of variable $v = (u\sqrt{\frac{\Delta}{N+1}} + y - x_0)/\sqrt{\Delta \frac{k}{N+1}}$, we can transform $J_4$ into the integral $\int_{\mathbb{R}} |v| e^{-v^2/2} dv$, which is finite and independent of $u$, $N$ and $k$. Combining these findings with (A.16), we can conclude that $|\frac{\partial F}{\partial u}(u, k, N)| \leq \frac{D}{\sqrt{N+1}}$ for some constant $D$ independent of $u$, $N$ and $k$.

By the fundamental theorem of calculus, we have $C(u, N) = C(0, N) + \int_0^u \frac{\partial C}{\partial u}(v, N) \, dv$. Using the bounds on $\frac{\partial F}{\partial u}$ we have just obtained, we find that the following inequalities hold:

$$|C(u, N) - C(0, N)| \leq \sum_{k=1}^{N-1} \frac{1}{\sqrt{k(N-k)}} \int_0^u \left| \frac{\partial F}{\partial u}(v, k, N) \right| dv$$

$$\leq \frac{Du}{\sqrt{N+1}} \sum_{k=1}^{N-1} \frac{1}{\sqrt{k(N-k)}}.$$

Since the sum on the right-hand side is convergent and since we know that the limit $L := \lim_{N \to \infty} C(0, N)$ exists, it follows from this inequality that $\lim_{N \to \infty} C(u, N) = L$ for every $u$ in $\mathbb{R}$. This proves (b) and (c).

Statement (d) follows from the dominated convergence theorem if we can show that the function $C(u, N)$ is bounded for $u \in \mathbb{R}$ and all integers $N$. This is a consequence of the definition of $C(u, N)$ given in (A.11), the fact that $F(u, k, N)$ is bounded [see (A.12)] and statement (a) which was proven above. This concludes the proof of the lemma in case (1), where $f_\delta(x, y) = h(x)((y - x)^2 - \delta)$, because, by (c) and (d), we get $\lim_{N \to \infty} B(u, N) = (\lim_{N \to \infty} C(u, N)) \int_{\mathbb{R}} (u^2 - 1) e^{-u^2/2} \, du = 0$, which is equivalent to the statement in (A.8).

Cases (2) and (3) are much simpler. Again, what we need to show is that (A.8) holds for the corresponding choices of $f_\delta(x, y)$. By introducing the substitution $u = \frac{x-y}{\sqrt{\delta}}$, where $\delta = \frac{\Delta}{N+1}$, the absolute value of the integral in (A.8) of case (2) [resp. (3)] is transformed to $\sqrt{\frac{\Delta}{N+1}} \int_{\mathbb{R}} |u|^3 e^{-u^2/2} C(u, N) \, du$ [resp., $\sqrt{\frac{\Delta}{N+1}} \int_{\mathbb{R}} |u| e^{-u^2/2} C(u, N) \, du$], where the function $C(u, N)$ is as in (A.11), while the function $F(u, k, N)$ is as in (A.12) with the exception of the integrand $h(u\sqrt{\frac{\Delta}{N+1}} + y, y)$, which now becomes a bounded function of two variables. It is clear that the equality in (A.8) will hold as soon as we see that the function $C(u, N)$ is bounded for all $u \in \mathbb{R}$ and all integers $N$. This follows from definition (A.11), from statement (a) on page 22 and from the inequality $|F(u, k, N)| \leq D_1 \int_{\mathbb{R}} \exp(-\frac{(x_\Delta - y)^2}{2\Delta}) \, dy$, which holds for some constant $D_1$. This concludes the proof of the lemma. $\quad\square$

## APPENDIX B: PROOF OF THEOREM 3

As mentioned in the discussion following Theorem 3, it follows form Lemma 1 in Beskos et al. (2006) and expressions (A.1) and (A.2) in Appendix A that the Radon–Nikodym derivative $d\mathbb{Q}_\theta^{\Delta,x,y}/d\mathbb{W}^{\Delta,x,y}$ exists and is of the required form. The bound in (a) follows from the Cauchy–Schwarz inequality on the Hilbert space $L^2(\mathbb{W}^{\Delta,x,y})$ since the Brownian bridge is a Gaussian process with bounded variance and the function $G$ has at most polynomial growth. In other words, the random variable $G(W_{t_1}, \ldots, W_{t_M})$ is an element in $L^2(\mathbb{W}^{\Delta,x,y})$, as is the Radon–Nikodym derivative since the drift $\mu(\cdot, \theta)$ is bounded above and below.

Part (b) of the theorem is a consequence of the inequality in part (a), Jensen's inequality and the elementary observation $0 \leq \mathbb{E}[X^2] - \mathbb{E}[X]^2 \leq \frac{(b-a)^2}{4}$, which holds for any random variable $X$ that takes values in a bounded interval $[a, b]$, applied to $X := \frac{L_\theta^\Delta}{\mathbb{E}_{\mathbb{W}^{\Delta,x,y}}[L_\theta^\Delta]} - 1$. Here, the variable $L_\theta^\Delta$ is as defined in Theorem 3. This concludes the proof.

**Acknowledgments.** We would like to thank Yacine Aït-Sahalia, Martin Crowder, Helmut Elsinger, Axel Gandy, Alois Geyer, Klaus Pötzelberger, Leopold Sögner, Osnat Stramer and Ruey Tsay for useful discussions, Nick Bingham for detailed comments on an earlier version of the manuscript and the anonymous referee for noting that Theorem 1 holds in much greater generality than initially stated.

DEPARTMENT OF MATHEMATICS
IMPERIAL COLLEGE LONDON
HUXLEY BUILDING
180 QUEEN'S GATE
LONDON SW7 2AZ
UNITED KINGDOM
E-MAIL: a.mijatovic@imperial.ac.uk

FINANCE GROUP
WARWICK BUSINESS SCHOOL
SCARMAN ROAD
COVENTRY CV4 7AL
UNITED KINGDOM
E-MAIL: Paul.Schneider@wbs.ac.uk